\documentclass[11pt]{amsart}
\usepackage{amsmath,amsfonts,amsthm,mathrsfs}
\usepackage{amssymb}

\usepackage[unicode,bookmarks]{hyperref}
\usepackage[dvipsnames]{xcolor}
\hypersetup{colorlinks=true,citecolor=NavyBlue,linkcolor=Brown,urlcolor=Orange}

\usepackage{tikz-cd,adjustbox}
\usepackage[arrow,curve,matrix]{xy}

\usepackage{booktabs}

\newif\ifdraft
\draftfalse

\counterwithin{equation}{subsection}
\counterwithout{subsection}{section}
\counterwithin{figure}{subsection}

\newtheorem{theorem}[equation]{Theorem}
\newtheorem*{theorem*}{Theorem}
\newtheorem{lemma}[equation]{Lemma}
\newtheorem*{lemma*}{Lemma}
\newtheorem{corollary}[equation]{Corollary}
\newtheorem{proposition}[equation]{Proposition}
\newtheorem*{proposition*}{Proposition}

\theoremstyle{definition}

\newtheorem*{definition*}{Definition}
\newtheorem{remark}[equation]{Remark}

\newtheorem{example}[equation]{Example}
\newtheorem*{example*}{Example}
\newtheorem*{problem*}{Problem}

\theoremstyle{plain}

\newcounter{intro}

\newtheorem{intro-conjecture}[intro]{Conjecture}
\newtheorem{intro-corollary}[intro]{Corollary}
\newtheorem{intro-theorem}[intro]{Theorem}
\newtheorem{intro-proposition}[intro]{Proposition}

\newcommand{\theoremref}[1]{\hyperref[#1]{Theorem~\ref*{#1}}}
\newcommand{\lemmaref}[1]{\hyperref[#1]{Lemma~\ref*{#1}}}
\newcommand{\definitionref}[1]{\hyperref[#1]{Definition~\ref*{#1}}}
\newcommand{\propositionref}[1]{\hyperref[#1]{Proposition~\ref*{#1}}}
\newcommand{\conjectureref}[1]{\hyperref[#1]{Conjecture~\ref*{#1}}}
\newcommand{\corollaryref}[1]{\hyperref[#1]{Corollary~\ref*{#1}}}
\newcommand{\exampleref}[1]{\hyperref[#1]{Example~\ref*{#1}}}

\newcommand{\Q}{\mathbb{Q}}
\newcommand{\C}{\mathbb{C}}

\renewcommand{\O}{\mathcal{O}}
\newcommand{\I}{\mathcal{I}}
\newcommand{\cP}{\mathcal{P}}
\newcommand{\F}{\mathcal{F}}

\newcommand{\ch}{{\rm ch}}
\DeclareMathOperator{\Pic}{Pic}

\DeclareMathOperator{\red}{red}

\DeclareMathOperator{\codim}{codim}

\makeatletter
\let\old@caption\caption
\renewcommand*{\caption}[1]{%
	\setcounter{figure}{\value{equation}}%
	\stepcounter{equation}%
	\old@caption{#1}\relax%
}
\makeatother

\begin{document}

\title{Generic vanishing and singularities of theta divisors}

\author{Fanjun Meng}
\address{Department of Mathematics, University of California San Diego, 9500 Gilman Drive \# 0112, La Jolla, 92093-0112, USA}
\email{f2meng@ucsd.edu}

\author{Mihnea Popa}
\address{Department of Mathematics, Harvard University, 1 Oxford Street, Cambridge, MA 02138, USA}
\email{mpopa@math.harvard.edu}

\thanks{FM was partially supported by the Simons Collaboration Grant on Moduli of Varieties. MP was partially supported by the NSF grant DMS-2401498.}

\begin{abstract}
We study the connection between the generic vanishing conjecture for subvarieties of principally polarized abelian varieties and various predictions about the singularities of theta divisors.
We deduce new cases of a conjecture on the multiplicities of points on irreducible theta divisors, as well as Beauville's conjecture characterizing Jacobians by the size of their theta singularity loci.
\end{abstract}

\maketitle

\makeatletter
\newcommand\@dotsep{4.5}
\def\@tocline#1#2#3#4#5#6#7{\relax
  \ifnum #1>\c@tocdepth 
  \else
    \par \addpenalty\@secpenalty\addvspace{#2}%
    \begingroup \hyphenpenalty\@M
    \@ifempty{#4}{%
      \@tempdima\csname r@tocindent\number#1\endcsname\relax
    }{%
      \@tempdima#4\relax
    }%
    \parindent\z@ \leftskip#3\relax
    \advance\leftskip\@tempdima\relax
    \rightskip\@pnumwidth plus1em \parfillskip-\@pnumwidth
    #5\leavevmode\hskip-\@tempdima #6\relax
    \leaders\hbox{$\m@th
      \mkern \@dotsep mu\hbox{.}\mkern \@dotsep mu$}\hfill
    \hbox to\@pnumwidth{\@tocpagenum{#7}}\par
    \nobreak
    \endgroup
  \fi}
\def\l@section{\@tocline{1}{0pt}{1pc}{}{\bfseries}}
\def\l@subsection{\@tocline{2}{0pt}{25pt}{5pc}{}}
\makeatother



\section{Introduction}

All throughout we work over $\C$. The main technical result of the paper is a stronger version of a conjecture made by Pareschi and the second author, about the equivalence of two notions of regularity for ideal sheaves on abelian varieties, designed according to \cite{PP1} to correspond to the notion of Castelnuovo-Mumford $2$-regularity in projective space, namely the implication (3) $\implies$ (2) in \cite[Conjecture A]{PP3}.\footnote{It was already shown there that (2) implies (3), and is equivalent to (4), a condition in terms of $M$-regularity.} Note that while geometric nondegeneracy was an assumption in the conjecture, here we see that it is not needed, and in fact the stronger version of Ran nondegeneracy follows from the hypothesis.

\begin{intro-theorem}\label{thm:main}
 Let $(A, \Theta)$ be an indecomposable principally polarized abelian variety, and $X \subset A$ a proper closed subscheme. If $\I_X (2\Theta)$ satisfies $IT_0$, then $\I_X (\Theta)$ is a $GV$-sheaf. Moreover, X is nondegenerate and integral.
\end{intro-theorem}

The case when $X$ is a divisor was obtained earlier in \cite{Meng}. For the technical notions in the theorem above, please see \S\ref{scn:vanishing}. 
If a coherent sheaf $\F$ on an abelian variety is $GV$, then $\F (\Theta)$ satisfies $IT_0$; see \cite[Proposition 3.1]{PP4}. The converse however is not true for arbitrary $\F$; see \cite[Example 3.2]{PP3}. Therefore the conjecture/theorem above is a property that is specific to twisted ideal sheaves. 
The (a priori stronger)  $GV$ property for $\I_X (\Theta)$ is more useful for applications, as it is related to a Schottky-type conjecture formulated in \cite{PP3} (namely (2) $\implies$ (5) in that same Conjecture A in \emph{loc. cit.}), which has been verified in a few cases in \cite{PP3}, \cite{CMPS}, \cite{Chen}.

\begin{intro-conjecture}\label{conj:GV}
 Let $(A, \Theta)$ be an indecomposable principally polarized abelian variety of dimension $g$. If  $X \subset A$ is a geometrically nondegenerate closed subscheme of pure dimension\footnote{Note that the geometric nondegeneracy and pure dimension assumptions can now be removed by Theorem \ref{thm:main}.} $1 \le d \le g-2$, then $\I_X (\Theta)$ is a $GV$-sheaf if and only if one of the following holds: 
 \begin{enumerate}
\item $(A, \Theta)$ is the Jacobian of a smooth curve $C$, and $X$ is $+$ or $-$ a translate of $W_d (C)$.
\smallskip
\item $(A, \Theta)$ is the intermediate Jacobian of a smooth cubic threefold $Y$, and $X$ is $+$ or $-$ a translate of the Fano surface $F$ of lines on $Y$.
\end{enumerate}
\end{intro-conjecture}

Although not directly relevant to this paper, it is worth mentioning that both items in this conjecture should also be equivalent to $X$ having minimal class, i.e. $[X] = \theta^{g-d}/ (g-d)!$, by a well-known conjecture of Debarre \cite{Debarre1}, which is the most significant conjecture in this direction. It is shown in \cite{PP3} that every $GV$-subscheme, i.e. one satisfying the conditions in Conjecture \ref{conj:GV}, has minimal class. Therefore the conjecture on minimal classes can be considered as input in the results below; note that besides the Matsusaka-Ran criterion, which implies that it holds for curves, it is also known to hold in dimension up to four, by \cite{Ran}. For more in this direction see e.g. the Introduction of \cite{CMPS}.

While results about Conjecture \ref{conj:GV} do rely specifically on the $GV$ property,  the $IT_0$ property for $I_X (2\Theta)$ is easier to verify in practice, and as a consequence,  Theorem \ref{thm:main} unlocks a number of significant  applications regarding the singularities of theta divisors.
The main one is related to the following well-known conjecture of the multiplicities of points on irreducible theta divisors. The first part has been folklore for some time, see for instance \cite[Conjecture 5.12]{Grushevsky}, while the stronger second statement is formulated as \cite[Question 4.7]{Casalaina1}.

\begin{intro-conjecture}\label{conj:multiplicity}
 Let $(A, \Theta)$ be an indecomposable principally polarized abelian variety of dimension $g$. Then:
\begin{enumerate}
\item for every $x \in \Theta$ we have ${\rm mult}_x (\Theta) \le \frac{g+1}{2}$.
\smallskip
\item more generally, we have
$$\dim \Sigma_m (\Theta) \le g - 2m +1,$$
where $\Sigma_m (\Theta) : = \{x \in \Theta ~|~ {\rm mult}_x (\Theta) \ge m\}$.
\end{enumerate}
\end{intro-conjecture}

\noindent
{\bf Remarks.\,} (1) (\emph{Equality case})
Note that (2) implies that in (1) we can also add that the upper bound ${\rm mult}_x (\Theta) = \frac{g+1}{2}$, if attained at all, can only be attained at finitely many points. In fact we can conjecture something more precise: the upper bound is attained at exactly one point, and only on the Jacobian of a hyperelliptic curve of odd genus, at the multiple of a $g_2^1$, or at the origin of the intermediate Jacobian of a smooth cubic threefold.

\noindent
(2) (\emph{Stronger Conjecture \ref{conj:multiplicity}(2)})
Schnell and Yang predict in \cite[Conjecture 1.11]{SY} a more precise statement in part (2) of the conjecture above, namely that for all $m\ge 2$ we have
\begin{equation}\label{eqn:SY}
\dim \Sigma_m (\Theta) \le g - 2m,
\end{equation}
unless $A$ is a hyperelliptic Jacobian, or the intermediate Jacobian of a smooth cubic threefold.

Conjecture \ref{conj:multiplicity} is known in dimension up to $5$ by Casalaina-Martin's work on singularities of theta divisors on Prym varieties \cite{Casalaina2}, and it is also known when $\Theta$ has isolated singularities by \cite[Theorem I]{MP}. Using a strategy outlined by the second author in a 2018 Oberwolfach lecture \cite{Popa}, based on the use of Hodge ideals as in \cite{MP}, and which was conditional on a result like Theorem \ref{thm:main}, we show: 

\begin{intro-theorem}\label{thm:conj-implication}
(1) Conjecture \ref{conj:GV} implies Conjecture \ref{conj:multiplicity}.  

\noindent
(2) Moreover, assuming $(A, \Theta)$ is not a hyperelliptic Jacobian, Conjecture \ref{conj:GV} implies (\ref{eqn:SY}), unless possibly when $g +1 = 2m$ and $\Sigma_m (\Theta)$ is a single reduced point, which is $2$-torsion for $\Theta$ symmetric.
\end{intro-theorem}

The exceptional condition in (2) is satisfied by the intermediate Jacobian of a smooth cubic threefold. Showing that it does not happen 
otherwise is equivalent to resolving the question in Remark (1). Other than this case, we see that Conjecture \ref{conj:GV} also implies the Schnell-Yang conjecture, and the fact that equality can happen in Conjecture \ref{conj:multiplicity}(1) at only at most one point, which is a $2$-torsion point if $\Theta$ is symmetric.

Moving on to unconditional statements, thanks to results in \cite{PP3}, \cite{CMPS} and \cite{Chen} we get:

\begin{intro-corollary}\label{cor:mult-unconditional}
Conjecture \ref{conj:multiplicity} holds for $g \le 6$, and for arbitrary $g$ when $\dim {\rm Sing}(\Theta) \le 2$ or $\dim {\rm Sing}(\Theta) = g-3$.
\end{intro-corollary}

Recall that an irreducible theta divisor is normal, hence in any case 
$$\dim {\rm Sing} (\Theta) \le g-3.$$
For a more precise statement about the unconditional cases, please see Corollary \ref{cor:dim6}. This also gives an alternative proof for the case $g \le 5$.

An upper bound on the dimension of $\Sigma_m (\Theta)$ that goes in the direction of Conjecture \ref{conj:multiplicity}(2) is given in Theorem \ref{thm:generalized-multiplicity} (which 
is independent of Theorem \ref{thm:main} or Conjecture \ref{conj:GV}) and Corollary \ref{cor:unconditional}. A special case of the latter that is easy to state is:

\begin{intro-theorem}\label{thm:m=3}
Conjecture \ref{conj:multiplicity}(2) holds in arbitrary dimension for $m \le 3$, or when $\Theta$ is $1$-Du Bois.\footnote{This means that its minimal exponent satisfies 
$\widetilde{\alpha} (\Theta) \ge 2$, or  equivalently its first Hodge ideal satisfies $I_1 (\Theta) = \O_A$.}
\end{intro-theorem}

The method we use for the proof of Theorem \ref{thm:conj-implication} also leads to a proof of the following Schottky-type conjecture originating with Beauville \cite{Beauville}; see also Debarre's \cite[Conjecture 3]{Debarre2} and the surrounding discussion. It was verified in  dimension at most $5$  in \cite[Proposition 6.6]{Beauville}, using the theory of generalized Prym varieties.

\begin{intro-theorem}\label{thm:beauville}
If a ppav $(A, \Theta)$ of dimension $g$ satisfies $\dim {\rm Sing} (\Theta) = g -3$, then $(A, \Theta)$ is the Jacobian of a hyperelliptic curve.
\end{intro-theorem}

It is of course known that all hyperelliptic Jacobians satisfy this property.
It is also not hard to show that if $\Theta$ is reducible, then $\dim {\rm Sing} (\Theta) = g-2$. (See \cite[Corollary 2]{EL} for a converse to this statement, and more.) Therefore the conjecture/theorem fully describes the case of maximal dimension for singular loci of irreducible theta divisors.
As shown in \cite[\S4.3]{Debarre2} (see also \cite{HJ} for more details), it also implies a conjecture of Debarre \cite[Conjecture 2]{Debarre2} on the base loci of polarizations.

\begin{intro-corollary}\label{cor:ACP}
If $(A, L)$ be a polarized abelian variety, then ${\rm Bs}(L)$ is generically reduced along any component of codimension $2$.
\end{intro-corollary}

A stronger version, unfortunately with an incorrect proof, was first formulated in \cite{ACP}, and used there to obtain an algebro-geometric proof of the characterization of Jacobians as those PPAVs whose theta function satisfies the KP-equation, as well as a characterization via degenerate Kummer trisecants.
See also \cite{Pareschi} for further comments and work in this direction. We are informed by E. Arbarello and G. Pareschi that  Corollary \ref{cor:ACP} suffices to fix the proofs of these applications.

For the applications of Theorem \ref{thm:main}, the main technical tool used in the paper is the zero-scheme of the first Hodge ideal $I_1 (\Theta)$ of \cite{MP}. The method applies, however, to the center of minimal exponent of \cite{SY} as well, answering a problem raised there (and shown in dimension one) in a few additional cases; see Corollary \ref{cor:center-min-exp} for the precise statement.

\medskip

\noindent
{\bf Acknowledgements.} 
The second author would like to thank G. Pareschi for many conversations over the years. His ideas have shaped many of the questions and techniques in this paper.  We thank him and E. Arbarello for pointing out a consequence to their work on Schottky-type problems, as well as S. Grushevsky and J. M\textsuperscript{c}Kernan for helpful discussions and comments. R. Yang has kindly shared a few comments and arguments that improved the picture in \S\ref{scn:center}. We also thank S. Huang and Z. Jiang, who have shared their work \cite{HJ} written at the same time as ours, where they also show Theorem \ref{thm:main},  and Beauville's conjecture Theorem \ref{thm:beauville}, using the same general strategy.

AI was used for only cleaning up and catching small errors. Note however that some of the unconditional statements here, as well as in \cite{HJ}, rely on the recent preprint \cite{Chen}, which provides an AI proof of Conjecture \ref{conj:GV} for $GV$-surfaces in ppavs of dimension $g \ge 6$.

\section{Preliminaries}

\subsection{Vanishing properties}\label{scn:vanishing}
Let $(A,\Theta)$ be a complex principally polarized abelian variety (ppav), and let $\widehat A=\Pic^0(A)$.
The cohomological support loci of a coherent sheaf $\F$ on $A$ are
\[
V^i(\F):=\{\alpha\in\widehat A \mid H^i(A,\F\otimes\alpha)\ne0\}.
\]
We say that $\F$ \emph{satisfies $IT_0$} if
$$
H^j(A,\mathcal \F\otimes\alpha)=0 \,\,\,\,{\rm for ~all} \,\,\,\,\alpha\in\widehat A \,\,\,\,{\rm and ~all}\,\,\,\,j>0,
$$
i.e. $V^i(\F)=\emptyset$ for all  $i>0$.  We say that $\F$ is a \emph{$GV$-sheaf} if 
$$\codim_{\widehat A}V^i(\F)\ge i \,\,\,\,\,\,{\rm for~all}\,\,\,\,i >0.$$ 
A closed subscheme $X\subseteq A$ is called a $GV$-subscheme if $\I_X(\Theta)$ is a $GV$-sheaf.

\subsection{Theta-dual}\label{scn:theta-dual}
For a reduced closed subscheme $Y \subseteq A$, we denote its theta-dual by $V(Y)$. Set-theoretically,
$$V (Y ) =  \{a \in A ~|~ Y \subseteq \Theta + a\}.$$
In particular $Y - V(Y) \subseteq \Theta$. For a detailed study of theta-duals, please see \cite{PP3}; in particular, they are endowed there 
with a scheme structure, using the Fourier-Mukai transform.

\subsection{Pontrjagin product}
Let $p_1,p_2:A\times A\to A$ be the natural projections, and let
$$
m\colon A \times A \to A, \,\,\,\,m(x,y)=x+y
$$
be the addition map. For coherent sheaves $\mathcal E,\mathcal G$ on $A$, denote $\mathcal E\boxtimes\mathcal G
  =p_1^*\mathcal E\otimes p_2^*\mathcal G$, and define
$$
\mathcal E*\mathcal G=m_*(\mathcal E\boxtimes\mathcal G)
$$
to be their \emph{Pontrjagin product}, with the derived version $\mathcal E*^{\mathbf R}\mathcal G$.
Thus the ordinary and derived Pontrjagin products agree whenever
$R^j m_*(\mathcal E\boxtimes\mathcal G)=0$ for $j>0$.
More about sheaf-theoretic Pontrjagin products can be found in \cite[\S5]{PP2}.

By interchanging the factors via $\tau(x,y)=(y,x)$ and the equality $m\circ\tau=m$, we get a natural isomorphism
$$
\varphi_{\mathcal E,\mathcal G}:\mathcal E*\mathcal G
  \xrightarrow{\sim}\mathcal G*\mathcal E.
$$
If $a:\mathcal E\to\mathcal E'$ and
$b:\mathcal G\to\mathcal G'$ are morphisms, it is immediate to check that
\begin{equation}\label{eqn:commutation}
\varphi_{\mathcal E',\mathcal G'}\circ(a*b)
  =(b*a)\circ \varphi_{\mathcal E,\mathcal G}.
\end{equation}

\begin{example}\label{ex:pontrjagin}
(1) Let $a \in A$ and $M = \O(\Theta_a)$. We note that 
$$
R^i m_*(M\boxtimes M)=0 \,\,\,\,\,\,{\rm for ~all}\,\,\,\, i >0,
$$
which implies that $M*M \simeq M*^{\mathbf R} M$. (This will be useful later on.)
Indeed, for  $t\in A$ we identify $m^{-1}(t)$ with $A$ by $x\mapsto(x,t-x)$, and
we write $\sigma_t(x)=t-x$. Under this identification we have
$$
(M\boxtimes M)|_{m^{-1}(t)}
  \simeq  M\otimes\sigma_t^*M.
$$
This is an ample line bundle on $A$ for each $t$, hence it has vanishing higher cohomology, and therefore the cohomology and base change theorem 
applies to imply the claim.

\noindent
(2) In addition to the notation in (1), let now $X$ be a closed subscheme of $A$, and let $\F = \I_X \otimes M$. Assume in addition that $\I_X (2\Theta)$ satisfies $IT_0$. We again have
$$
R^i m_*(\F\boxtimes M)=0 \,\,\,\,\,\,{\rm for ~all}\,\,\,\, i >0,
$$
and consequently $\F*M \simeq \F*^{\mathbf R} M$. This time we have 
$$
(\F\boxtimes M)|_{m^{-1}(t)}
  \simeq  \F\otimes\sigma_t^*M.
$$
Since this is of the form $\I_X (2 \Theta_b)$ for some $b$, it has vanishing higher cohomology by our assumption, and therefore the same argument applies.
\end{example}

\subsection{Hodge ideals}\label{scn:Hodge}
We continue to consider a complex ppav $(A, \Theta)$. Recall from \cite{MP} that there is a Hodge filtration $F_\bullet \O_A (*\Theta)$ on the $\mathcal{D}_A$-module 
$\O_A(*\Theta)$ of meromorphic functions on $A$ with poles of arbitrary order. In each degree it is contained in the pole order filtration 
$P_k \O_A(*\Theta) = \O_A \big( (k+1)\Theta\big)$, hence for each $k \ge 0$ we have 
$$F_k \O_A (*\Theta) = \O_A((k+1)\Theta) \otimes I_k (\Theta),$$
where $I_k (\Theta)$ is the $k$-th Hodge ideal of $\Theta$. For $k = 0$ we have 
$$I_0 (\Theta) \simeq \I \big((1-\epsilon) \Theta\big)$$
the multiplier ideal of the $\Q$-divisor $(1-\epsilon) \Theta$, where $0 < \epsilon \ll 1$ is a rational number; see \cite[Proposition 10.1]{MP}. 
This holds for any reduced divisor on a smooth variety, but in our case  a well known result of Koll\'ar \cite[Theorem 17.13]{Kollar} says then that $I_0 (\Theta) \simeq \O_A$, i.e. the pair $(A, \Theta)$ is log canonical.

 In this paper we will deal with the Hodge ideal $I_1 (\Theta)$. We will need:
  
 \begin{theorem}\label{thm:Hodge}
 We have:
 \begin{enumerate}
 \item For every irreducible closed subset $W\subset A$ of codimension $r$, if $m= {\rm mult}_W (\Theta)$, then $I_1 (\Theta) \subseteq \I_W^{(q)}$, where
 $q= {\rm min}~\{ m-1, 2m-r \}$. (Here $ \I_W^{(q)}$ is the $q$-th symbolic power of the ideal sheaf $\I_W$.)
 \item $H^i \big(A, I_1 (\Theta)\otimes \O_A (2 \Theta) \otimes \alpha\big) = 0$ for all $i > 0$ and all $\alpha \in \widehat{A}$, i.e. $I_1 (\Theta)\otimes \O_A(2 \Theta)$
 satisfies $IT_0$.
 \end{enumerate}
 
 \end{theorem}

The first is a special case of \cite[Theorem E]{MP}, a local property that holds on any smooth variety, while the second is a special case of \cite[Theorem 28.2]{MP}, which is in turn the better behaved abelian varieties version of the Hodge ideal analogue of Nadel's vanishing theorem.

\section{GV-subschemes}

\subsection{Vanishing away from theta-dual}
The first thing to show is that all the cohomological support loci of $\I_X (\Theta)$ are subsets of the theta-dual of $X_{\rm red}$ (see \S\ref{scn:theta-dual}), under the $IT_0$ assumption. (Eventually we will see that $X$ must in fact be integral.)

\begin{proposition}\label{prop:containment}
Let $(A, \Theta)$ be a ppav, and let $X$ be a closed subscheme of $A$. If $\I_X (2 \Theta)$ satisfies $IT_0$, then we have
$$H^i (A, \I_X (\Theta_a)) = 0 \,\,\,\,{\rm for~all}\,\,\,\,a \not\in V(X_{\rm red}) \,\,\,\,{\rm and~all}\,\,\,\,i \ge 0.$$
In other words, for all $i \ge 0$ we have 
$$V^i \big (\I_X (\Theta)\big) \subseteq V (X_{\rm red}).$$
\end{proposition}
\begin{proof}
Let $a \not\in V(X_{\rm red})$, and set  
$$M= \O_A (\Theta_a) \,\,\,\,\,\,{\rm and}\,\,\,\,\,\, \F=\mathcal I_X\otimes M.$$ 
We are aiming to show that 
$$H^i (A, \F) = 0  \,\,\,\,{\rm for~all}\,\,\,\,i \ge 0.$$

By the choice of $a$, there exist a closed point $p\in X$, and a section $s \in H^0 (A, M)$, such that $s(p)\ne 0$. By abuse of notation, we denote by $M(p)$ the fiber of $M$ at $p$, thought of as a skyscraper sheaf at $p$ or a vector space over $\C$ depending on the specific setting. Consider
$$
f \colon \mathcal F\hookrightarrow M,
\qquad
g\colon M\rightarrow M(p),
$$
i.e. the natural inclusion and evaluation maps. Clearly  $g\circ f=0$.

We consider the composition of Pontrjagin products
\begin{equation}\label{eqn:composition}
\F*M \overset{f*{\rm id}}{\longrightarrow} M*M \overset{{\rm id} *g}{\longrightarrow} M*M(p),
\end{equation}
\emph{Claim:} this composition is equal to $0$.

Let us first conclude the proof assuming the Claim. It is immediate to check that the composed morphism above also factors as 
$$\F*M \overset{{\rm id} *g}{\longrightarrow} \F*M(p) \overset{f*{\rm id}}{\longrightarrow} M*M(p).$$
Here the second map is injective; it suffices to have the injectivity of the map $\F \boxtimes M(p) \to M \boxtimes M(p)$ on $A\times A$, but this can be identified on $A \times \{p\}$ with the inclusion $\F\hookrightarrow M$ tensored over $\C$ with a vector space, which stays injective.
Therefore we deduce that 
$${\rm id} *g \colon \F*M \longrightarrow \F*M(p)$$
is the $0$ map. On the other hand, by Example \ref{ex:pontrjagin}(2), for all $i$ we have 
$$
H^i(A,\mathcal F*M) \simeq H^i(A\times A,\mathcal F\boxtimes M) \simeq H^i(A,\mathcal F)\otimes H^0(A,M),
$$
as the higher cohomology of $M$ vanishes. Similarly, 
$$
H^i(A,\mathcal F*M(p))
  \simeq H^i(A\times A,\mathcal F\boxtimes M(p))
  \simeq H^i(A,\mathcal F)\otimes M(p).
$$
Moreover, the map on cohomology induced by ${\rm id}*g$ is identified with
$$
H^i(A,\mathcal F)\otimes H^0(A,M)
  \xrightarrow{\;1\otimes {\rm ev}_p\;}
H^i(A,\mathcal F)\otimes M(p),
$$
where ${\rm ev}_p$ is the evaluation of global sections at $p$. Now $H^0(A,M)=\mathbb C s$ and $M(p)$ are both one-dimensional, 
and ${\rm ev}_p(s)=s(p)\ne 0$, so ${\rm ev}_p$ is an isomorphism. It follows that $H^i(A,\mathcal F)=0$ for all $i\geq 0$, as required.

It remains to prove the Claim, i.e. the vanishing of the composition in ($\ref{eqn:composition}$). To this end, we first show that the involution 
$\varphi_{M, M}$ on $M * M$ induced by the exchange map $(x,y) \to (y, x)$ on $A \times A$ is equal to the identity. To see this, we first note that 
$E = M*M$ is a simple vector bundle on $A$. Indeed, by Example \ref{ex:pontrjagin}(1) we have $E\simeq M*^{\mathbf R} M$, and by \cite[\S3, formula (3.7)]{Mukai} we know that the Fourier-Mukai transform 
$$
\Phi  \colon D^b(A)\longrightarrow D^b(\widehat A),
$$
interchanges Pontrjagin and tensor product. Therefore 
$$
\Phi_A(E)\simeq\widehat M^{\otimes 2}.
$$
Here $\widehat{M}$ is the total Fourier-Mukai transform of $M$, which is well known to be a line bundle on $\widehat{A}$ (the line bundle associated to a translate of $-\Theta$). The endomorphism group of a line bundle is $\C$, and since $\Phi$ is an equivalence \cite[Theorem 2.2]{Mukai}, we obtain 
${\rm End}_A (E) \simeq \C$.

Since $E$ is simple, we have that $\varphi_{M, M} = \lambda \cdot {\rm id}_E$ for some $\lambda \in \C$. Moreover, $\varphi_{M, M}$ fixes the nontrivial 
section $s \boxtimes s$ of $E$, since this is preserved by the exchange involution. Hence $\lambda = 1$.

We use this to deduce our claim, namely that the morphism 
$$\Psi = ({\rm id}*g)\circ(f* {\rm id})\colon \F*M\to  M*M(p)
$$
is $0$. By applying the commutation relation ($\ref{eqn:commutation}$) to  ${\rm id}_M$ and $g$, we get
$$
\varphi_{M,M(p)}\circ({\rm id}* g)=(g * {\rm id})\circ \varphi_{M,M}.
$$
Using $\varphi_{M,M}= {\rm id}_E$, we finally get
$$
\varphi_{M,M(p)}\circ \Psi =(g * {\rm id}_M)\circ(f * {\rm id}_M)  =(g\circ f)* {\rm id}_M =0.
$$
But $\varphi_{M,M(p)}$ is an isomorphism, so $\Psi =0$.
\end{proof}

\begin{corollary}\label{cor:GV-via-nondeg}
Under the hypothesis of Proposition \ref{prop:containment}, if $X$ is geometrically nondegenerate, then $I_X (\Theta)$ is a $GV$-sheaf.
\end{corollary}
\begin{proof}
Note that for $i > 0$ and any $a \in A$ we have 
$$H^{i+1}  (A, \I_X (\Theta_a)) \simeq H^i   (X, \O_X (\Theta_a)),$$
which vanishes for $i \ge d +1$, where $d = \dim X$. Therefore we only need to concern ourselves with $V^i \big (\I_X (\Theta)\big)$ for $i \le d+1$, and show that their dimension is at most $g - i$.

Now since $X$ is geometrically nondegenerate, so is $X_{\rm red}$ and by \cite[Lemma 4.4]{PP3} we have 
$$\dim V (X_{\rm red}) \le g - d- 1.$$
Combined with Proposition \ref{prop:containment}, this immediately implies what we want.
\end{proof}

Geometric nondegeneracy will be treated in the next section, but it is automatic in the $0$-dimensional case, and this already has an interesting consequence.

\begin{corollary}\label{lem:MP}
Let $Y \subseteq A$ be a $0$-dimensional closed subscheme such that $\I_Y (2 \Theta)$ satisfies $IT_0$. Then $Y$ is exactly one reduced point.
\end{corollary}
\begin{proof}
According to Corollary \ref{cor:GV-via-nondeg}, $\I_X (\Theta)$ is a $GV$-sheaf. Thus for general $a \in A$ we have $H^1 (A, \I_X (\Theta_a)) = 0$, and therefore the restriction map
$$H^0 (A, \O_A (\Theta_a)) \to H^0 (X, \O_X (\Theta_a))$$
is surjective. This implies that the length of $X$ is $1$. 
\end{proof}

This provides in particular a different approach to \cite[Theorem 29.2]{MP}  in \emph{loc. cit.}, which is a proof of the multiplicity Conjecture \ref{conj:multiplicity} when $\Theta$ has isolated singularities. A careful reading of the proof in \emph{loc. cit.} shows that it relies precisely on the statement of Corollary \ref{lem:MP}, proved there 
using Theorem \ref{thm:Hodge}, plus the fact that the linear series $|2 \Theta|$ separates neither points nor tangent vectors on $A$.

\subsection{Geometric nondegeneracy}
The second ingredient for the proof of Theorem \ref{thm:main} is a result that combined with the work of Schreieder \cite{Schreieder} ensures the nondegeneracy of a closed subscheme $X$ for which $\I_X (2\Theta)$ satisfies the $IT_0$ property.

\begin{proposition}\label{prop:complementary-theta-summands}
Let $(A,\Theta)$ be an indecomposable principally polarized abelian variety of dimension $g$, and $X\subset A$ a proper closed subscheme of dimension $d$.  Assume that $\I_X(2\Theta)$ satisfies $IT_0$. Then for every integral closed subscheme $Y$ of $X_{\red}$ of dimension $d$, there exists an integral closed subscheme $Z\subseteq V(X_{\red})$ of dimension $g-d-1$ such that 
\[
Y-Z=\Theta.
\]
In particular, $Y$ is (Ran) nondegenerate.
\end{proposition}

\begin{proof}
The last assertion follows for $Y$ from the first. If  $1 \le d \le g-2$, it follows by \cite[Theorem 1.1(1)]{Schreieder}. If  $d = g-1$, we get that $Y$ is a translate of $\Theta$, in which case it is clear. If $d = 0$,  nondegeneracy is automatic.

We now focus on the first assertion. Let $\cP$ be the normalized Poincar\'e bundle on $A\times\widehat A$. Under the isomorphism $A\simeq\widehat A$ induced by the principal polarization, we can regard $\cP$ as a line bundle on $A\times A$. 

Set 
$$\theta =c_1(\O_A(\Theta))\in H^2(A,\Q) \,\,\,\,{\rm  and} \,\,\,\, \ell =c_1(\cP)\in H^2(A\times A,\Q).$$
Let 
$$\delta:A\times A\to A, \,\,\,\,(a,b)\mapsto a-b$$
be the difference map. 
 By the construction of $\cP$ in \cite[p.78]{Mumford}, our choice of the isomorphism $A\simeq\widehat A$, and the fact that $\cP$ is trivial when restricted to 
 $A\times\{0\}$ or $\{0\}\times A$, we have that
\[
\ell=\delta^*\theta-p_1^*\theta-p_2^*\theta \in H^1(A,\Q)\otimes_{\Q}H^1(A,\Q)\subseteq H^2(A\times A,\Q).
\]

Denote the two natural projections by $p_1,p_2 \colon A\times A\to A$, and consider the Fourier-Mukai transform of $\I_X (\Theta)$, namely
\[
C:= \Phi \bigl(\I_X (\Theta)\bigr) = Rp_{2*}\bigl(p_1^*\I_X(\Theta)\otimes\cP\bigr).
\]
By Proposition \ref{prop:containment}, we have $C |_{A \smallsetminus V(X_{\red})} = 0$.

We now proceed to some Riemann-Roch calculations. Recall that there is an induced Fourier-Mukai transform at the level of cohomology with rational coefficients, defined on 
$H^\bullet (A, \Q)$ by 
$$
\mathscr F(\alpha)=p_{2*}\bigl(p_1^* \alpha \cup e^\ell\bigr).
$$
Using the fact that $\theta^g = g!$, we have
$$
\mathscr F(e^\theta)
  =p_{2*}\bigl(e^{\delta^*\theta-p_2^*\theta}\bigr)
  =e^{-\theta}.
$$
Using the triviality of the tangent bundle of $A$, as well as the standard exact sequence involving $\I_X$ and $\O_X$, Grothendieck-Riemann-Roch implies
$$
\ch(C)=e^{-\theta}
  -\mathscr F\bigl(\ch(\mathcal O_X(\Theta))\bigr).
$$
Now since the support of $\mathcal O_X(\Theta)$ has dimension $d$, we have
$\ch_r(\mathcal O_X(\Theta))=0$ for $r<g-d$. Moreover, since the bidegree of $\ell$ is $(1,1)$ and the reverse Fourier-Mukai transform is the inverse of $\mathscr F$ up to sign and the standard involution, for all $j$ we have
$$
\mathscr F\bigl(H^{2j}(A,\Q)\bigr)
  =  H^{2(g-j)}(A,\Q).
$$
Thus for every $j>d$ we deduce that
\[
\ch_{j}(C)=(-1)^{j}\frac{\theta^{j}}{j!}, 
\]
and in particular $\ch_{d+1}(C)=(-1)^{d+1}\frac{\theta^{d+1}}{(d+1)!}.$

Since $C |_{A \smallsetminus V(X_{\red})} = 0$, the Chern character $\ch (C)$ in the Chow ring of $A$ restricts to 
$0$ on this open set, and in particular $\ch_{d+1}(C)$ is represented by a rational cycle class with support in $V(X_{\red})$. Passing to cohomology, we obtain an expression
$$
\frac{\theta^{d+1}}{(d+1)!}=\sum_j r_j[Z_j],
$$
where $r_j\in\Q$ and each $Z_j\subseteq V(X_{\red})$ is a closed subvariety of dimension $g-d-1$; here $[W]$ denotes the cohomology class of a closed subvariety $W$.

Now fix a $d$-dimensional component $Y\subseteq X_{\mathrm{red}}$ as in the hypothesis, and set
$$
\gamma
  =\delta_*\left(p_1^* [Y]\cup p_2^*\frac{\theta^{d+1}}{(d+1)!}\right)
  =\sum_j r_j\,\delta_*[Y\times Z_j]
  \in H^2(A,\mathbb Q).
$$
For the Pontrjagin product
$u* v=m_*(p_1^*u\cup p_2^*v)$ of cohomology classes, we have the standard identity
$$
\frac{\theta^{d+1}}{(d+1)!} *
\frac{\theta^{g-1}}{(g-1)!}
  =(g-d)\frac{\theta^d}{d!}.
$$
See for instance \cite[Corollary 16.5.8]{BL}.
 The automorphism $f\colon A\times A\to A\times A$ defined by $(a,b)\mapsto(b,a-b)$ satisfies $m\circ f=p_1$, $p_1\circ f=p_2$, and $p_2\circ f=\delta$. Thus we deduce that
\begin{align*}
\int_A\gamma\cup\frac{\theta^{g-1}}{(g-1)!}&=\int_{A\times A}p_1^*[Y]\cup p_2^*\frac{\theta^{d+1}}{(d+1)!}\cup\delta^*\frac{\theta^{g-1}}{(g-1)!}\\
&=\int_{A\times A}m^*[Y]\cup p_1^*\frac{\theta^{d+1}}{(d+1)!}\cup p_2^*\frac{\theta^{g-1}}{(g-1)!}\\
&=\int_A[Y]\cup\left(\frac{\theta^{d+1}}{(d+1)!} * \frac{\theta^{g-1}}{(g-1)!}\right)=\frac{g-d}{d!}\int_Y \theta^d>0.
\end{align*}
Thus we deduce that $\delta_*[Y\times Z_j]\ne0$ for some $j$ and therefore
\[
\dim(Y-Z_j)=\dim(Y\times Z_j)=g-1.
\]
On the other
hand,
$$
Y-Z_j\subseteq X_{\mathrm{red}}-V(X_{\rm red})\subseteq\Theta.
$$
Since the theta divisor of an indecomposable ppav is integral, we deduce that $Y-Z_j=\Theta$. Taking  $Z=Z_j$ concludes the proof.
\end{proof}

\subsection{Proof of Theorem \ref{thm:main}}
We assume that $X$ is a proper closed subscheme of an indecomposable ppav $(A, \Theta)$ such that $\I_X (2 \Theta)$ satisfies $IT_0$. 
By Proposition \ref{prop:complementary-theta-summands} we know that any component $Y$ of $X_{\rm red}$ with $\dim Y = \dim X$ is nondegenerate, 
hence also geometrically nondegenerate. We deduce that  $X_{\rm red}$ is also geometrically nondegenerate by \cite[Lemma 4]{Schreieder}. Hence so is $X$, because this property only depends on the reduced scheme structure.  Then by Corollary \ref{cor:GV-via-nondeg} we deduce that $\I_X (\Theta)$ is a $GV$-sheaf. As $X$ is a geometrically nondegenerate closed $GV$-subscheme of $A$, \cite[Theorem 2(1)]{Schreieder} implies that $X$ is integral. Finally, this upgrades $X$ to being nondegenerate, since it is equal to $Y$.

\begin{remark}
Note that once it is shown that $X$ is a $GV$-scheme, and therefore automatically nondegenerate by Theorem \ref{thm:main}, then 
various other nice properties of $X$ follow. We have already seen that $X$ is integral by \cite{Schreieder}. It is also known that $X$ is normal, and even has rational singularities when $\dim X = 2$, by \cite{JP}. Further generation properties relying only on the $GV$-condition were studied in \cite{LT}.
\end{remark}

\section{Singularities of theta divisors}

\subsection{Multiplicities of points, and multiplicity loci}
Let $(A, \Theta)$ be a complex ppav.  We will make use of the Hodge ideal $I_1 (\Theta)$, see \S\ref{scn:Hodge}. From here on, write 
$$I_1 (\Theta) = \I_Z,$$
where $Z$ is the closed subscheme of $A$ determined by this ideal.

We start with the following result towards Conjecture \ref{conj:multiplicity}, which is independent of the generic vanishing conjecture or Theorem \ref{thm:main}.
Note first that when $\Theta$ is irreducible, it has rational singularities by \cite[Theorem 1]{EL}, and in particular it is normal. Therefore we have
$$\dim Z \le \dim {\rm Sing} (\Theta) \le \dim \Theta - 2 = g -3.$$
For this reason,  Conjecture \ref{conj:multiplicity}  is immediate for $m =2$.

\begin{theorem}\label{thm:generalized-multiplicity}
If $(A, \Theta)$ is an indecomposable ppav of dimension $g$, then for all $m \ge 2$ we have
$$\dim \Sigma_m (\Theta) \le {\rm max}~\{g - 2m + 1, \dim {\rm Sing}(Z) \},$$
with the convention that $\dim \emptyset = - \infty$.\footnote{Note that this result only uses statements from \cite{MP}, and nothing about $GV$-schemes or minimal classes.}
\end{theorem}
\begin{proof}
For $m =2$ the result is clear by the discussion above, so we may assume $m \ge 3$.
Recall that we consider $\Sigma_m (\Theta)$ with its reduced scheme structure, and let $W$ be an irreducible component. 
Denote 
$$r = \codim_A W \,\,\,\,\,{\rm and} \,\,\,\,\,\, e = {\rm mult}_W (\Theta).$$
In particular we have $e \ge m \ge 3$.
If $r \ge 2m - 1$, or equivalently $\dim W \le g -2m +1$, we are in good shape already. So we may assume $r \le 2m - 2$.

With this data, we apply Theorem \ref{thm:Hodge}(1) to obtain
$$\I_Z \subseteq \I_W^{(q)}, \,\,\,\,\,{\rm for}\,\,\,\,q = {\rm min}~\{e -1, 2e -r\},$$
where $\I_W^{(q)}$ is the $q$-th symbolic power of $\I_W$. Given our numerical assumptions, we have 
$q \ge 2$, so that at the very least 
$$\I_Z \subseteq \I_W^{(2)}.$$
Let $x \in W$ be a general point, at which $W$ is smooth. Since the symbolic and usual powers agree at smooth points, we have
$$\I_{Z, x} \subseteq \I_{W, x}^2 \subseteq \frak{m}_x^2,$$
where the latter is the maximal ideal of $x$ in $A$. This implies that $x$ is a singular point of $Z$, and therefore $W \subseteq {\rm Sing}(Z)$.
\end{proof}

Studying the case $\dim Z = 0$ relies on Corollary \ref{lem:MP}, which says that $Z$ is a reduced point. 
Next we begin to make use of Theorem \ref{thm:main}, if $\dim Z \ge 1$. The sheaf $\I_Z (2\Theta)$ satisfies $IT_0$  by Theorem \ref{thm:Hodge}(2), and consequently  $Z$ is a $GV$-scheme by Theorem \ref{thm:main}. It follows that  $Z$ is itself normal by \cite[Lemma 7.8]{JP}, and so overall
$$\dim {\rm Sing}(Z)  \le \dim Z - 2 \le g -5.$$

Putting everything together, we get a conclusion that includes Theorem \ref{thm:m=3}.

\begin{corollary}\label{cor:unconditional}
Conjecture \ref{conj:multiplicity}(2) holds in arbitrary dimension for 
$$m \le \frac{g +1 - \dim {\rm Sing}(Z)}{2}.$$
Hence it always holds for $m \le 3$. It holds for all $m$ if $Z$ is smooth or empty.
\end{corollary}

The condition $Z = \emptyset$ is equivalent to having $I_1 (\Theta) = \O_A$, which is in turn equivalent to the condition $\widetilde{\alpha} (\Theta) \ge 2$ for the minimal exponent of $\Theta$ by \cite{Saito}. This is known as the $1$-Du Bois condition; equivalently, the pair $(A, \Theta)$ is $1$-log canonical in the terminology of \cite{MP}. This is stronger than rationality; for instance, in the case of Jacobians it holds for Brill-Noether general curves, by \cite[Theorem 5.1]{SY}.

\medskip

We next prove the main result in this direction.

\begin{proof}[Proof of Theorem \ref{thm:conj-implication}]
The proofs of (1) and (2) could be packaged together a bit more efficiently. We keep them separate however, in order to emphasize exactly what is needed 
for each. 

\noindent
(1)  We need to show that the generic vanishing Conjecture \ref{conj:GV} implies the strong multiplicity Conjecture \ref{conj:multiplicity}.

If $Z$ is smooth or empty, the conjecture holds by Corollary \ref{cor:unconditional}. In particular this applies when $\dim Z = 0$, as we may apply Corollary \ref{lem:MP} to deduce that $Z$ is smooth.
Assume now that $\dim Z > 0$. Then by Theorem \ref{thm:main} we have that $I_Z (\Theta)$ is a $GV$-sheaf, and $Z$ is integral and geometrically nondegenerate. Conjecture \ref{conj:GV}  then leads to two possibilities. 

On one hand, $(A, \Theta)$ could be the intermediate Jacobian of a smooth cubic threefold; this however cannot be the case, since its theta divisor has only one singular point, 
contradicting the assumption on $\dim Z$. On the other hand, $(A, \Theta)$ could be the Jacobian of a smooth projective curve of genus $g$. Then the bound in Conjecture \ref{conj:multiplicity}(2) holds by Martens' theorem.

\noindent
(2) We may choose $\Theta$ to be symmetric. In this case $Z$ is symmetric as well, i.e. $Z = -Z$. Indeed, Hodge ideals commute with pullback by  automorphisms by \cite[Proposition 15.1]{MP}, and therefore in this case $I_1 (\Theta)$ is invariant under the involution $a \mapsto -a$.

We next claim that if $W \subseteq \Sigma_m (\Theta)$ is an irreducible component such that $\dim W \ge g - 2m +1$, then $W \subseteq Z$. 
We make use of the argument in the proof of Theorem \ref{thm:generalized-multiplicity}. Say 
$$r = \codim_A W \,\,\,\,\,{\rm and} \,\,\,\,\,\, e = {\rm mult}_W (\Theta).$$
We have $e \ge m$ and $2m -1 \ge r$.  Thanks to Theorem \ref{thm:Hodge}(1), we have
$$\I_Z \subseteq \I_W^{(q)}, \,\,\,\,\,{\rm for}\,\,\,\,q = {\rm min}~\{e -1, 2e -r\} \ge 1.$$
Hence $W \subseteq Z$. In particular this shows that we may not have $Z = \emptyset$.

If $d=\dim Z > 0$, then as in (1) $A$ must be the Jacobian of a curve $C$, and $Z$ a copy of $W_d (C)$ (up to sign and translation). It is however well known that if the property $Z = -Z$ holds for some $W_d (C)$, then $C$ must be hyperelliptic; see also \cite[Lemma 5.5]{Debarre1} for a more general statement.
This is excluded by our hypothesis.

Thus the only possibility that is left is $g +1 = 2m$ and $\dim Z = 0$. The latter implies that $Z$ is a single reduced point by Corollary \ref{lem:MP}. 
It is also a $2$-torsion point by its symmetry. Finally, by (1) and the discussion above, each component of $\Sigma_m (\Theta)$ is $0$-dimensional and contained in $Z$, hence $\Sigma_m (\Theta) = Z$.
\end{proof}

The following is a slightly stronger version of Corollary \ref{cor:mult-unconditional}.

\begin{corollary}\label{cor:dim6}
Conjecture \ref{conj:multiplicity} holds unconditionally when:
\begin{enumerate}
\item $g=\dim A \le 6$
\item More generally, when $\dim Z \le 2$ or $\dim Z \ge g-3$, in any dimension $g$; this holds in particular when $\dim {\rm Sing} (\Theta) \le 2$ or  $\dim {\rm Sing} (\Theta) \ge g-3$.
\end{enumerate}
\end{corollary}
\begin{proof}
It suffices to show (2). The case $\dim Z = 0$ is dealt with by Corollary \ref{lem:MP}. 
Next note that the case $\dim Z \le 2$ is technically equivalent to the case $\dim Z \ge g-3$. This is because by theta-duality, see 
\cite[Theorem 5.2]{PP3}, if $Z$ is a $GV$-scheme of dimension $d$, its theta-dual $W$ is a $GV$-scheme of dimension $g - d -1$, and the argument in the proof of 
Theorem \ref{thm:conj-implication} works for either. 

Hence all we need to know is that the $GV$-conjecture holds for all positive dimensional closed subschemes $X$ of dimension at most two. If $\dim X = 1$, this is known by \cite[Theorem C]{PP3}. If $\dim X = 2$, it is known if $g \le 5$ by \cite[Theorem A]{CMPS}, and if $g \ge 6$ by \cite[Corollary 1.3]{Chen}.

Finally, if $\dim {\rm Sing} (\Theta) \le 2$, then clearly $\dim Z \le 2$. On the other hand, if  $\dim {\rm Sing} (\Theta) \ge g-3$, then in fact 
$\dim {\rm Sing} (\Theta) = g-3$ since $\Theta$ is normal, in which case we also have $\dim Z = g -3$ by Lemma \ref{lem:g-3} below.
\end{proof}

\begin{remark}
The first case when we do not have the multiplicity conjecture unconditionally is when $g = 7$ and $\dim Z = \dim {\rm Sing}(\Theta) = 3$. This is the only case left to deal with in dimension $7$, and the conjecture is that it should not be possible (as it is not hard to see that it cannot happen on Jacobians, with $Z = W_3$ up to sign and translation).
\end{remark}

\begin{remark}[Equality case]
This refers to the Remark right after Conjecture \ref{conj:multiplicity} in the Introduction, regarding the maximal multiplicity case
$${\rm mult}_x (\Theta) = \frac{g+1}{2}.$$
On the Jacobian of a smooth curve it is well known that the multiplicity of each point on 
$\Theta$ is at most $(g+1)/2$, by the Riemann Singularity Theorem combined with Clifford's Theorem. Moreover, Clifford also tells us that the equality case 
happens only on hyperelliptic Jacobians, at a unique point (the multiple of the $g_2^1$) when $g$ is odd, and not at all when $g$ is even. It is also well known that on the intermediate Jacobian of a smooth cubic threefold the only singular point is $0 \in \Theta$, which also satisfies the equality case ${\rm mult}_0 (\Theta) = 3$. It it natural to conjecture that these are the only situations when equality holds.
\end{remark}

\subsection{Characterization of hyperelliptic Jacobians}
The methods discussed above also lead to a proof of Beauville's conjecture characterizing hyperelliptic Jacobians; see Theorem \ref{thm:beauville} and the preceding paragraph.
We start with a technical result:

\begin{lemma}\label{lem:g-3}
If $\dim {\rm Sing} (\Theta) = g-3$, then $\dim Z = g-3$. Moreover, ${\rm Sing}(\Theta)$ has a unique irreducible component of maximal dimension $g-3$, equal to $Z$. 
\end{lemma}
\begin{proof}
Clearly we have $Z \subseteq {\rm Sing}(\Theta)$. By  Theorem \ref{thm:main}, we also know that $Z$ is integral. It suffices therefore to show that if $W$ is a $(g-3)$-dimensional 
component of ${\rm Sing}(\Theta)$, then $W \subseteq Z$.

Denote by $m$ the generic multiplicity of $\Theta$ along $W$. Since the pair $(A, \Theta)$ is log canonical (see \S\ref{scn:Hodge}), it is standard that, for any $k \ge 2$, each irreducible component of $\Sigma_k (\Theta)$ has codimension at least $k$ in $A$; moreover, by \cite[Corollary 2]{EL}, the existence of a component of codimension precisely $k$ implies that $\Theta$ is reducible. Therefore in our case we have the inequality
$$g -3 = \dim W \le g - m -1,$$
and therefore $m = 2$. We now use the relationship between Hodge ideals and multiplicity loci  in Theorem \ref{thm:Hodge}(1) ; in our case, with codimension $r = 3$ and multiplicity $m =2$,  it implies that 
$$I_1 (\Theta) \subseteq I_W,$$
which is precisely what we want.
\end{proof}

\begin{proof}[Proof of Theorem \ref{thm:beauville}]
It is a standard fact that if $\Theta$ is reducible, then $\dim {\rm Sing} (\Theta) = g-2$. Therefore under our assumption $\Theta$ is irreducible.

As discussed in the Introduction, Beauville \cite{Beauville} proved this result for $g \le 5$, hence we may assume $g \ge 6$.
By Lemma \ref{lem:g-3}, we see that if $Z$ is the zero-scheme of the Hodge ideal $I_1 (\Theta)$, then $Z$ itself has dimension $g-3$. Moreover, as before 
Theorem \ref{thm:main} implies that $\I_Z (\Theta)$ is a $GV$-sheaf. Using  the main result of \cite{Chen}, the $GV$-conjecture holds for $Z$, and therefore $(A, \Theta)$ is a Jacobian. But on the Jacobian of a curve $C$ it is well known that $\dim {\rm Sing} (\Theta) = g-3$ if and only if $C$ is hyperelliptic.
\end{proof}

\subsection{Centers of minimal exponent}\label{scn:center}
Another closed subvariety that is amenable to the same techniques as above is the center of minimal exponent $Y$ of $\Theta$, defined by Schnell-Yang \cite[\S4.2]{SY}. 
From the point of view of this paper, $Y$ plays a role similar to that of the zero-scheme $Z$ of the ideal $I_1 (\Theta)$, but it is not hard to check that in general they are different. For instance by construction we always have $Y \neq \emptyset$ when $\Theta$ is singular.
This is also the place to note that some of the arguments in this paper can be equally applied to the zero-schemes of the higher multiplier ideals of \emph{loc. cit.}; they are indeed considered in \cite{HJ}.

It is shown in \cite[Theorem 1.12]{SY} that if $\dim Y =1$, then $Y$ must be a hyperelliptic curve, and if moreover there exists a component of the maximal allowed dimension $g -2m +1$ in some $\Sigma_m (\Theta)$, then $A$ is the Jacobian of $Y$.  Please see the interesting discussion in \S1.6 of \emph{loc. cit.} for the significance of such a statement.  Using the discussion in the paragraph above, we have:

\begin{corollary}\label{cor:center-min-exp}
Assume that  $\dim Y \ge 1$, and the minimal exponent of $\Theta$ satisfies $1 < \widetilde{\alpha} (\Theta) < 2$. If Conjecture \ref{conj:GV} holds, then $(A, \Theta)$ is the Jacobian of a hyperelliptic curve.\footnote{It is shown in \cite{SY} that every Jacobian satisfies $1 < \widetilde{\alpha} (\Theta) \le 2$, while for Brill-Noether general curves we have 
$\widetilde{\alpha} (\Theta) = 2$. Note also that the condition $1 < \widetilde{\alpha} (\Theta)$ is automatic, since it is equivalent to rational singularities.}

In particular, this holds unconditionally if $1 \le \dim Y \le 2$, or if $\dim Y = g-3$.
\end{corollary}
\begin{proof}
In the proof of \cite[Lemma 5.5]{SY} it is shown that under this minimal exponent hypothesis $\I_Y (2 \Theta)$ satisfies $IT_0$, so the methods we developed here apply: by Theorem \ref{thm:main}, $\I_Y (\Theta)$ is $GV$, and $Y$ is a geometrically nondegenerate (integral, normal) closed subvariety of $A$.

We obviously have $Y \subseteq {\rm Sing} (\Theta)$, and therefore $1\le \dim Y \le g -3$. Assuming Conjecture \ref{conj:GV}, we deduce that since $Y$ is positive dimensional, $(A, \Theta)$ is the Jacobian of a curve $C$. The conjecture also says that $Y$ must be $+$ or $-$ a translate of $W_d (C)$ where $d=\dim Y$. But we may also assume that $Y = -Y$, by choosing 
$\Theta$ to be symmetric, as the construction of the center of minimal exponent is invariant under the involution by negation. 
As in the proof of Theorem \ref{thm:conj-implication}(2), this implies that $C$ must be hyperelliptic.

The unconditional statement refers to the cases of the $GV$-conjecture we know to hold, as everywhere else in the paper. 
\end{proof}

Restricting to Jacobians $A = J(C)$, it is shown in \cite{SY} that if $C$ is hyperelliptic (when $\widetilde{\alpha} (\Theta) = 3/2$) or Brill-Noether general (when 
$\widetilde{\alpha} (\Theta) = 2$, the maximal possible value for Jacobians), then $Y = {\rm Sing} (\Theta)$. We can now also deduce what happens for non-hyperelliptic curves with $\widetilde{\alpha} (\Theta) < 2$; we thank Ruijie Yang for pointing this out, together with the argument and the subsequent example.

\begin{corollary}
If $A$ is the Jacobian of a non-hyperelliptic curve $C$, and $\widetilde{\alpha} (\Theta) < 2$, then $Y$ is a single reduced point, which is $2$-torsion when $\Theta$ is symmetric.
\end{corollary}
\begin{proof}
As in the previous proof, $\I_Y (\Theta)$ is $GV$, and $Y$ is integral and geometrically nondegenerate. By \cite[Theorem B]{PP3}, $Y$ is a subvariety of minimal class in $A$. 
Since $A$ is a Jacobian, by \cite[Theorem 5.1]{Debarre1}, $Y$ is a $W_d$ (up to sign and translation). But we also know that if $\Theta$ is chosen to be symmetric, then $Y = -Y$, hence if it were positive dimensional as before $C$ would be hyperelliptic. We deduce that $\dim Y = 0$, and $Y$ is a single reduced point, which must be $2$-torsion in the symmetric case.
\end{proof}

\begin{example}
It can also happen however to have $\widetilde{\alpha} (\Theta)  = 2$ for curves that are non-hyperelliptic but not Brill-Noether general. Indeed, if a non-hyperelliptic curve with $g \ge 4$ satisfies the property that  $\Theta$ has no singular $2$-torsion point, then $\widetilde{\alpha} (\Theta)  = 2$. This follows from the argument above, since $Y \subseteq {\rm Sing}(\Theta)$. This holds for instance for a trigonal genus $5$ curve obtained as the normalization of a plane quintic with a unique singularity which is an ordinary node, by \cite[Ch.V, p.207--208]{ACGH}.
\end{example}


\begin{thebibliography}{9}

\bibitem[ACP]{ACP}
E.~ Arbarello, G.~Codogni, G.~Pareschi, \emph{Characterizing Jacobians via the KP equation and via flexes
and degenerate trisecants to the Kummer variety: an algebro-geometric approach} J. reine angew. Math. \textbf{777} (2021), 251--271.

\bibitem[ACGH]{ACGH}
E.~Arbarello, M.~Cornalba, P. A.~Griffiths, and J.~Harris, \emph{Geometry of algebraic curves} Vol. I,
volume 267 of \emph{Grundlehren der mathematischen Wissenschaften [Fundamental Principles of Math-
ematical Sciences]}, Springer-Verlag, New York, 1985.

\bibitem[Be]{Beauville}
A.~Beauville, \emph{Prym varieties and the Schottky problem},
Invent. Math. \textbf{41} (1977), 149--196.

\bibitem[BL]{BL}
Ch.~Birkenhake and H.~Lange, \emph{Complex abelian varieties}, Second, Grundlehren der mathema-
tischen Wissenschaften [Fundamental Principles of Mathematical Sciences], vol. 302, Springer-Verlag, Berlin, 2004.

\bibitem[Ca1]{Casalaina1}
S.~Casalaina-Martin, \emph{Singularities of theta divisors in algebraic geometry}, In Curves and abelian
varieties, Contemp. Math. \textbf{465}, 25--43. Amer. Math. Soc., Providence, RI, 2008.

\bibitem[Ca2]{Casalaina2}
S.~Casalaina-Martin, \emph{Singularities of the Prym theta divisor}, Ann. of Math. (2) \textbf{170} (2009), no.1, 162--204.

\bibitem[CMPS]{CMPS}
S.~Casalaina-Martin, M.~Popa, and S.~Schreieder, \emph{Generic vanishing and minimal cohomology classes on abelian fivefolds}, J. Algebraic Geom. \textbf{27} (2018),  553--581.

\bibitem[Ch]{Chen}
Y.~Chen, \emph{Generic vanishing subschemes of codimension three}, preprint arXiv:2609.15391.

\bibitem[De1]{Debarre1}
O.~Debarre, \emph{Minimal cohomology classes and Jacobians}, J. Algebraic Geom. \textbf{4} (1995),  321--335.

\bibitem[De2]{Debarre2}
O.~Debarre, \emph{Update on: Singularities of divisors on abelian varieties}, lecture notes at https://perso.imj-prg.fr/olivier-debarre/exposes-talks/.

\bibitem[EL]{EL}
L.~Ein and R. Lazarsfeld, \emph{Singularities of theta divisors and the birational geometry of irregular varieties}, 
J. Amer. Math. Soc. \textbf{10} (1997), no.1,  243--258.

\bibitem[Gr]{Grushevsky}
S.~Grushevsky, \emph{The Schottky problem}, In \emph{Current developments in algebraic geometry}, Math. Sci. Res. Inst. Publ. \textbf{59},  
129--164. Cambridge Univ. Press, Cambridge, 2012.

\bibitem[HJ]{HJ}
S.~Huang and Z.~Jiang, \emph{$GV$-subschemes and hyperelliptic Jacobians}, preprint. 

\bibitem[JP]{JP}
Z.~Jiang and G.~Pareschi, \emph{Cohomological rank functions on abelian varieties}, Ann. Sci. \'Ec. Norm. Sup\'er. (4) \textbf{53} (2020), no. 4, 815--846.

\bibitem[Ko]{Kollar}
J.~Koll\'ar, \emph{Shafarevich maps and automorphic forms}, M.B. Porter Lectures, Princeton Univ. Press, 1995.

\bibitem[LT]{LT}
L.~Lombardi and S.~Tirabassi, \emph{$GV$-subschemes and their embeddings in principally polarized abelian varieties}, Math. Nachr. \textbf{288} (2015), no. 11-12, 1405--1412.

\bibitem[Me]{Meng}
F.~Meng, \emph{Generic vanishing conjecture for divisors}, preprint arXiv:2609.11377.

\bibitem[Muk]{Mukai}
S.~Mukai, \emph{Duality between $D(X)$ and $D(\widehat X)$ with its application
to Picard sheaves}, Nagoya Math. J. \textbf{81} (1981), 153--175.

\bibitem[Mu]{Mumford}
D.~Mumford, \emph{Abelian varieties}, Tata Institute of Fundamental Research Studies in Mathematics,
vol. 5, Tata Institute of Fundamental Research, Bombay; by Hindustan Book Agency, New
Delhi, 2008. With appendices by C. P. Ramanujam and Yuri Manin, Corrected reprint of the second (1974) edition.

\bibitem[MP]{MP}
M.~Musta\c t\u a and M.~Popa, \emph{Hodge ideals},  Mem. Amer. Math. Soc. \textbf{262} (2019), no.1268.

\bibitem[Pa]{Pareschi}
G.~Pareschi, \emph{Singularities of base loci on abelian varieties}, Atti Accad. Naz. Lincei Rend. Lincei Mat. Appl. \textbf{36} (2025), no. 2, 365--375.

\bibitem[PP1]{PP1}
G.~Pareschi and M.~Popa, \emph{Regularity on abelian varieties I}, J. Amer. Math. Soc. \textbf{16} (2003), 285--302.

\bibitem[PP2]{PP2}
G.~Pareschi and M.~Popa, \emph{Regularity on abelian varieties II: Basic results on linear series
and defining equations}, J. Algebraic Geom. \textbf{13} (2004), no.~1, 167--193.

\bibitem[PP3]{PP3}
G.~Pareschi and M.~Popa, \emph{Generic vanishing and minimal cohomology classes on abelian varieties},
Math. Ann. \textbf{340} (2008), 209--222.

\bibitem[PP4]{PP4}
G.~Pareschi and M.~Popa, \emph{Regularity on abelian varieties III: relationship with Generic Vanishing and applications}, 
in \emph{Grassmannians, Moduli Spaces and Vector Bundles} 2008, Clay Mathematics Proceedings \textbf{14}, Amer. Math. Soc., Providence RI, 2011, 141--167.

\bibitem[Po]{Popa}
M.~Popa, \emph{Connections between some conjectures on subvarieties of abelian varieties}, Oberwolfach Report  33/2018, 1993--1995.

\bibitem[Ra]{Ran}
 Z.~Ran, \emph{On subvarieties of abelian varieties}, Invent. Math. \textbf{62} (1981), 459--479.

\bibitem[Sa]{Saito}
M.~Saito, \emph{Hodge ideals and microlocal $V$-filtration}, preprint arXiv:1612.08667. 
	
\bibitem[Sch]{Schreieder}
S.~Schreieder, \emph{Decomposable theta divisors and generic vanishing}, Int. Math. Res. Not. IMRN \textbf{16} (2017), 4984--5009.

\bibitem[SY]{SY}
C.~Schnell and R.~Yang, \emph{Higher multiplier ideals}, J. reine angew. Math. \textbf{832} (2026), 153--197.

\end{thebibliography}
\end{document}